\documentstyle{article}

\newtheorem{r}{Remark}
\newtheorem{lm}{Lemma}

\begin{document}

\title
{\bf Kauffman Monoids}
 
\author{{\sc Mirjana Borisavljevi\' c}
\\University of Belgrade
\\Faculty of Transport and Traffic Engineering
\\Vojvode Stepe 305, 11000 Belgrade, Yugoslavia
\\email: mirjanab@afrodita.rcub.bg.ac.yu
\\[0.2cm]
{\sc Kosta Do\v sen} 
\\Mathematical Institute
\\Knez Mihailova 35, P.O. Box 367
\\11001 Belgrade, Yugoslavia, and
\vspace {0.1cm}
\\IRIT, University of Toulouse III  
\\31062 Toulouse cedex, France
\\email: kosta@mi.sanu.ac.yu
\\[0.2cm]
{\sc Zoran Petri\' c}
\\Mathematical Institute
\\Knez Mihailova 35, P.O. Box 367
\\11001 Belgrade, Yugoslavia
\\email: zpetric@mi.sanu.ac.yu}

\date{ }
\maketitle
\begin{abstract}
\noindent This paper gives a self-contained and complete proof of the
isomorphism of freely generated monoids extracted from Temperley-Lieb
algebras with monoids made of Kauffman's diagrams.
\\[0.1cm]
\end{abstract}

\vspace{0.4cm}
\setcounter {section}{-1}
\section{Introduction}
\noindent Kauffman monoids (i.e. semigroups with unit)
are multiplicative structures extracted
from Temperley-Lieb algebras, which play a very prominent role in
knot theory,
in low-dimensional
topology, in topological quantum field theories, in quantum groups and in
statistical mechanics. The label ``Temperley-Lieb" is derived from a
paper in this last
field  \cite{TL}.
The axioms of Kauffman monoids (see Section 1 below) are,
however, not stated explicitly in that paper. They may be found in
Jones' paper \cite [$\S$4.1.4]{J1}, and they
appear in full light in many works of
Kauffman's \cite [$\S$4]{K1}, \cite [pp. 431-433]{K2},
\cite [IV]{K3}, \cite [I.7, I.16, II.8]{K4a}, \cite [IV]{K5},
\cite [II]{K6},
\cite [$\S$7.2, pp. 8-9]{K7}, \cite [Section 3]{CKS}, \cite
[Section 6]{K8}---in particular,
in \cite {K4}, which is most detailed. Although it is tempting to call
the monoids in question ``Temperley-Lieb monoids", we
believe ``Kauffman monoids"
is a fairer denomination, doing
justice to the author who unearthed these structures
(cf. \cite [p. 324]{J2}).
                                                                   
Kauffman monoids are seldom, if ever, separated from
Temperley-Lieb algebras, which are much richer structures, and  
not much attention is paid usually to the completeness of the standard
axiomatics of Kauffman monoids with respect to their standard
geometrical interpretation. Although this completeness, which
consists in showing that the freely generated Kauffman monoids are
isomorphic to the connection monoids of \cite{K4}, is well known, we
have been unable to find in the literature a self-contained and
exhaustive proof. Among the papers we know,
Kauffman's \cite{K4} and \cite[Section 6]{K8},
together with an aside in Jones' paper \cite[$\S$4.1.4]{J1}, whose
content may also be found in \cite[$\S$2.8]{GHJ},
are nearest to the mark.
There are, however, 
some uncertainties in \cite{K4} concerning the appropriate normal
form for terms of Kauffman monoids (see the proof of Theorem 4.3
on p. 442,
where the normal form mentioned is not Jones', and is not unique; see also
p. 434 of the same paper).
All these papers still leave pretty much
work to the reader.

We have been readers who have done the work left to us.
The completeness proof for
Kauffman monoids we
are going to present here is, we believe, self-contained and
thorough.
Moreover, it presents some aspects of Kauffman monoids that may be novel
and worth recording. We have not encountered elsewhere the block
formulation of
freely generated Kauffman monoids of Section 1. We use this formulation to
reduce terms syntactically to normal form, which is an essential step
in the completeness proof. In this formulation of Kauffman monoids
all axioms are tied to reduction to the normal form. We
have also not found elsewhere a rigorous proof, such as we present in
Section 3,
that all elements of connection monoids can be generated from
what Kauffman called
{\it hooks} (and which we shall call {\it diapsides}).  The best
we have found is
a proof in \cite[Section 6]{K8}, which leaves some details to the reader,
and a sketchy proof in \cite[Theorem 26.10,
$\S$26, Chapter VIII, pp. 251-253]{PS}. In \cite[Proposition 4.1.3]{BJ}
one may find a proof of something more general, and somewhat more
complicated.
Our proof exhibits some difficulties,
which we think cannot be evaded,
and it also sheds some light on the normal form of terms. 
We have not found elsewhere Lemma 4, the key lemma of Section 4
and of the whole completeness proof.

\section{The monoids ${\cal K}_{n}$}
\noindent The {\it Kauffman monoid}  ${\cal K}_{n}$, for $n \ge 2$, has for
every $i \in \bf N$ such that $1 \le i \le n-1$ a generator $h^i$,
called a {\it diapsis}
(plural {\it diapsides}), and also the generator $c$, called the {\it circle}.
(Kauffman in \cite{K4} called diapsides ``hooks", and that is where the
label $h$ comes
from; our Greek neologism means ``double arc".) The terms of
${\cal K}_{n}$
are defined inductively by stipulating that the generators and {\bf 1}
are terms,
and that if $t$ and $u$ are terms, then $(tu)$ is a term (as usual, we
shall omit the outermost parentheses of terms).

The monoid ${\cal K}_{n}$ is freely generated from the generators above
so that the following equations hold between terms of ${\cal K}_{n}$:
\begin{tabbing}
\hspace{10em} \= $(1)$ \hspace{2em} \= ${\bf 1} t = t {\bf 1} = t$,
\\[0.1cm]
\> $(2)$ \> $t(uv)=(tu)v$,
\\[0.1cm]
\> $(h1)$ \> $h^i h^j = h^j h^i, \;\; \mbox{ for } |i-j| \ge 2,$
\\[0.1cm]
\> $(h2)$ \> $h^i h^{i \pm 1}h^i =  h^i,$
\\[0.1cm]
\>$(hc1)$ \> $h^i c = c h^i,$
\\[0.1cm]
\>$(hc2)$ \> $h^i h^i = c h^i.$
\end{tabbing}

For $1\le j \le i \le n-1$, let the {\it block} $h^{[i,j]}$ be defined as
$h^i h^{i-1}  \ldots  h^{j+1} h^j$. The block $h^{[i,i]}$, which is
defined as
$h^i$, will be called {\it singular}. Let $c^1$ be $c$, and let
$c^{l+1}$ be $c^{l}c$.

A term is in {\it Jones normal form} iff it is either of the form 
$c^l h^{[b_1,a_1]} \ldots h^{[b_k,a_k]}$ for $l,k \ge 1$, 
$a_1< \ldots <a_k$ and $b_1< \ldots <b_k$,
or of the form
$h^{[b_1,a_1]} \ldots h^{[b_k,a_k]}$ for $k \ge 1$, 
$a_1< \ldots <a_k$ and $b_1< \ldots <b_k$,
or of the form $c^l$ for $l \ge 1$, or it is the term {\bf 1}
(see \cite [$ \S$4.1.4]{J1}). For the sake of definiteness, we require
that in the Jones normal form all parentheses are associated to the left
(but another arrangement of parentheses would do as well). In the reduction
to Jones normal form of Lemma 1 below we shall not bother with trivial
considerations concerning parentheses. The associativity equation $(2)$
guarantees that we can move them at will (as it dispensed us from
writing parentheses in $(h2)$).

That every term of ${\cal K}_{n}$ is equal to a term in Jones normal form
will be demonstrated with the help of an alternative formulation of
${\cal K}_{n}$,
called the {\it block formulation}, which is obtained as follows. Besides
the circle $c$, we take as generators the blocks $h^{[i,j]}$
instead of the diapsides, we generate terms with these generators,
{\bf 1} and multiplication, and to the equations
$(1)$ and $(2)$ we add the equations
      
\begin{tabbing}
\hspace{5em} \= $(h \rm I)$ \hspace{2em} \= $\mbox {for} \; \; j \ge k+2,
\; \; 
h^{[i,j]}h^{[k,l]} = h^{[k,l]}h^{[i,j]},$
\\[0.2cm]
\> $(h \rm II)$ \> $\mbox {for} \; \; i \ge l \; \mbox {and} \;
|k-j|=1, \; \; h^{[i,j]}h^{[k,l]} = h^{[i,l]},$
\\[0.2cm]
\> $(hc \rm I)$ \> $h^{[i,j]} c = c h^{[i,j]},$
\\[0.2cm]
\> $(hc \rm II)$ \> $h^{[i,j]}h^{[j,l]}= c h^{[i,l]}.$
\end{tabbing}

We verify first that with $h^i$ defined as the singular block $h^{[i,i]}$
the equations
$(h1),(h2), (hc1)$ and $(hc2)$ are instances of the new equations of
the block formulation. The equation
$(h1)$ is $(h \rm I)$ for $i=j$ and $k=l$, the equation $(h2)$
is $(h \rm II)$ for
$i=j=l$ and $k=i+1$, or $i=k=l$ and $j=i-1$, the equation $(hc1)$
is $(hc \rm I)$
for $i=j$, and the equation $(hc2)$ is $(hc \rm II)$ for $i=j=l$.
We also have to verify
that in the new axiomatization we can deduce the definition of $h^{[i,j]}$  
with diapsides replaced by singular blocks; namely, we have to verify 
$$h^{[i,j]}=h^{[i,i]}h^{[i-1,i-1]} \ldots h^{[j+1,j+1]}h^{[j,j]},$$
which readily follows from $(h \rm II)$ for $j=k+1$. To
finish showing that the block formulation of ${\cal K}_{n}$
is equivalent to
the old formulation, we have to verify that with blocks
defined via diapsides we can deduce
$(h {\rm I}), (h {\rm II}), (hc {\rm I})$ and $(hc \rm II)$
from the old equations,
which is a straightforward exercise.

We can deduce the following equations in ${\cal K}_{n}$ for 
$j+2 \le k$:
\begin{tabbing}
\hspace{3em} \= $(h \rm III.1)$ \hspace{2em} \=
$h^{[i,j]}h^{[k,l]} = h^{[k-2,l]}h^{[i,j+2]} \; \; 
\mbox {if} \; \;
i \ge k \; \; \mbox {and} \; \; j \ge l,$
\\[0.2cm]
\> $(h \rm III.2)$ \> $h^{[i,j]}h^{[k,l]} = h^{[i,l]}h^{[k,j+2]}
\; \; \; \; \; \; \mbox {if} \; \;
i < k \; \; \mbox {and} \; \; j \ge l,$
\\[0.2cm]
\> $(h \rm III.3)$ \> $h^{[i,j]}h^{[k,l]} = h^{[k-2,j]}h^{[i,l]} \;
\; \; \; \; \; \mbox {if} \; \;
i \ge k \; \; \mbox {and} \; \; j < l,$
\end{tabbing}

\noindent which is also pretty straightforward. Then we 
prove the following lemma.

\begin{lm}
Every term of ${\cal K}_{n}$ is equal in ${\cal K}_{n}$ to a term
in Jones normal form.
\end{lm}                             
{\it Proof.} We shall give a reduction procedure that transforms
every term into a term in Jones normal form, every reduction step
being justified by an equation of ${\cal K}_{n}$. (In logical jargon, we
establish that this procedure is strongly normalizing---namely, that
any sequence of reduction steps terminates in a term in normal form.)

Take a term in the block formulation of ${\cal K}_{n}$, and let
subterms of this term of the forms

\[\begin{array}{l}
h^{[i,j]} h^{[k,l]}, \; \; \mbox {for} \; \; i \ge k \; \; \mbox {or} \; \;
j \ge l,
\\[0.1cm] h^{[i,j]}c,
\\[0.1cm] {\bf 1}t, t{\bf 1}
\end{array}\]

\noindent be called {\it redexes}. A reduction of the first sort
consists in replacing
a redex of the first form by the corresponding term on the right-hand
side of one
of the equations $(h {\rm I}),(h {\rm II}),(hc \rm II),
(h {\rm III.1}),(h \rm III.2)$
and
$(h \rm III.3)$. (Note that the terms on the left-hand sides of these
equations cover
all possible redexes of the first form, and the conditions of these
equations
exclude each other.) A reduction of the second sort consists in replacing
a redex of the second form by
the right-hand side of $(hc \rm I)$, and, finally, a reduction of the
third sort consists in replacing
a redex of one the forms in the third line by $t$, according to
equation $(1)$.

Let the {\it weight} of a block $h^{[i,j]}$ be $i-j+2$. For any
subterm  $h^{[i,j]}$ of
a term $t$ in the block formulation of ${\cal K}_{n}$, let $\rho (h^{[i,j]})$
be the number of subterms $h^{[k,l]}$ of $t$ on the right-hand side
of $h^{[i,j]}$ such that $i \ge k$ or $j \ge l$. The subterms $h^{[k,l]}$
are not necessarily immediately on the right-hand side of $h^{[i,j]}$ as in
redexes of the first form:
they may also be separated by other terms. For any subterm $c$ of a term $t$
in the block formulation ${\cal K}_{n}$, let $\tau (c)$ be the number of
blocks on the left-hand side of this $c$.

The {\it complexity measure} of a term $t$ in the block formulation is 
$\mu (t)= (n_1,n_2)$ where $n_1 \ge 0$ is the
sum of the weights of all the blocks in $t$,
and $n_2 \ge 0$
is the sum of all the numbers $\rho (h^{[i,j]})$ for all blocks
$h^{[i,j]}$ in $t$ plus the sum of all the
numbers $\tau (c)$ for all circles $c$ in $t$ and plus the number of
occurrences of {\bf 1} in $t$. The ordered pairs $(n_1,n_2)$
are well-ordered lexicographically.

Then we check that if $t'$ is obtained from $t$ by a reduction,
then $\mu (t')$
is strictly smaller than $\mu (t)$. With reductions of the first sort
we have that
if they are based on $(h \rm I)$, then $n_2$ diminishes while
$n_1$ doesn't change, and if they are based on
the remaining equations, then $n_1$ diminishes. 
With reductions of the
second and third sort, $n_2$ diminishes while $n_1$ doesn't change.

So, by induction on the complexity measure, we obtain that every term 
is equal to a term without redexes, and it is easy to see that a term
is without redexes iff it is in Jones normal form. {\it q.e.d.}
\vspace{0.2cm}

Note that for a term  $c^l h^{[b_1,a_1]} \ldots h^{[b_k,a_k]}$
in Jones normal form the number $a_i$ is strictly smaller than all indices
of diapsides on the right-hand side of $h^{a_i}$, and $b_i$
is strictly greater than all indices
of diapsides on the left-hand side of $h^{b_i}$. The following remark is an
immediate consequence of that.

\begin{r}
If in a term in Jones normal form a diapsis $h^i$ occurs more than once,
then in between any
two occurrences of $h^i$ we have an occurrence of $h^{i+1}$ and an
occurrence of $h^{i-1}$.
\end{r}

A normal form dual to Jones' is obtained with blocks $h_{[i,j]}$ where
$i \le j$,
which are defined as $h^i h^{i+1} \ldots h^{j-1} h^j$. Then in
$c^l h_{[a_1,b_1]} \ldots h_{[a_k,b_k]}$
we require that $a_1> \ldots >a_k$ and $b_1> \ldots >b_k$.
The length of this new normal form will be the same as the length of Jones'.
As a matter of fact, we could take as a term in normal form many other terms
of the same reduced length as terms in the Jones normal form. For all these
alternative normal forms we can establish the property of Remark 1.

\section{The {\it n}-diagrams}
\noindent A one-manifold with boundary is a topological space
whose points have open
neighbourhoods homeomorphic to the real intervals $(-1,1)$ or $[0,1)$,
the boundary points
having the latter kind of neighbourhoods. For $n \ge 2$ a natural number and
$a>0$ a real number, let $R_{n,a}$ be the rectangle $[0, n+1] \times [0,a]$.
Let $[0, n+1] \times \{ a \} \subseteq R_{n,a}$ be the {\it top} of
$R_{n,a}$
and let $[0, n+1] \times \{ 0 \} \subseteq R_{n,a}$ be the {\it bottom} of
$R_{n,a}$.

An {\it n-diagram} $D$ is a compact one-manifold with boundary with a
finite number
of connected components embedded in a rectangle $R_{n,a}$ such that the
intersection of $D$ with the top of $R_{n,a}$ is
$t(D)= \{ (i,a)\;|\; i \in {\bf N} \cap [1,n] \}$, the
intersection of $D$ with the bottom of $R_{n,a}$ is
$b(D)= \{ (i,0)\;|\; i \in {\bf N} \cap [1,n] \}$,
and $t(D) \cup b(D)$ is the set of boundary points of $D$.

It follows from this definition that every {\it n}-diagram has
$n$ components homeomorphic to $[0,1]$, which are called
{\it threads}, and a finite number of components homeomorphic
to $S^1$, which are called {\it circular components}. The threads and 
the circular
components make all the connected components of an {\it n}-diagram.
All these components are mutually disjoint.

Every thread has two {\it end points} that belong to the boundary
$t(D) \cup b(D)$. When one of these end points is in $t(D)$
and the other in $b(D)$, the thread is {\it transversal}. 
A transversal thread is {\it vertical} when the
first coordinates of its end points are equal.
A thread that is not transversal is a {\it cup} when both of 
its end points are in $t(D)$,
and it is a {\it cap} when they are both in $b(D)$.
It is clear that the following holds.

\begin{r}
In every {\it n}-diagram the number of cups is equal to the number of caps. 
\end{r}

For example, an 11-diagram in $R_{11,10}$ looks as follows:

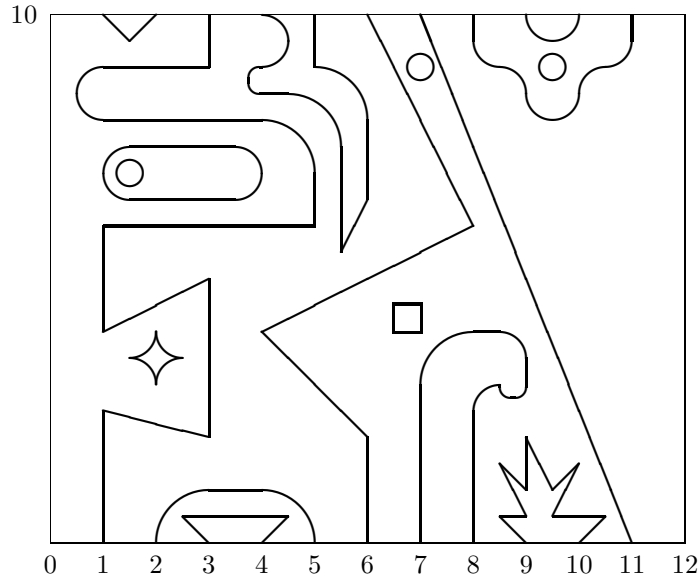
\begin{center}
\begin{picture}(240,240)
{\linethickness{0.05pt}
\put(0,20){\line(1,0){240}}
\put(0,20){\line(0,1){200}}
\put(0,220){\line(1,0){240}}
\put(240,20){\line(0,1){200}}}

\put(-5,220){\makebox(0,0)[r]{$10$}}
\put(0,15){\makebox(0,0)[t]{$0$}}
\put(20,15){\makebox(0,0)[t]{$1$}}
\put(40,15){\makebox(0,0)[t]{$2$}}
\put(60,15){\makebox(0,0)[t]{$3$}}
\put(80,15){\makebox(0,0)[t]{$4$}}
\put(100,15){\makebox(0,0)[t]{$5$}}
\put(120,15){\makebox(0,0)[t]{$6$}}
\put(140,15){\makebox(0,0)[t]{$7$}}
\put(160,15){\makebox(0,0)[t]{$8$}}
\put(180,15){\makebox(0,0)[t]{$9$}}
\put(200,15){\makebox(0,0)[t]{$10$}}
\put(220,15){\makebox(0,0)[t]{$11$}}
\put(240,15){\makebox(0,0)[t]{$12$}}

\thicklines
\put(20,20){\line(0,1){50}}
\put(20,70){\line(4,-1){40}}
\put(60,60){\line(0,1){60}}
\put(20,100){\line(2,1){40}}
\put(20,100){\line(0,1){40}}
\put(20,140){\line(1,0){80}}
\put(20,200){\line(1,0){40}}
\put(60,200){\line(0,1){20}}
\put(20,140){\oval(160,80)[tr]}
\put(20,190){\oval(20,20)[l]}

\put(20,220){\line(1,-1){10}}
\put(40,220){\line(-1,-1){10}}

\put(50,160){\oval(60,20)}

\put(30,160){\circle{10}}

\put(50,100){\oval(20,20)[bl]}
\put(30,100){\oval(20,20)[br]}
\put(50,80){\oval(20,20)[tl]}
\put(30,80){\oval(20,20)[tr]}

\put(70,20){\oval(60,40)[t]}

\put(60,20){\line(-1,1){10}}
\put(50,30){\line(1,0){40}}
\put(80,20){\line(1,1){10}}

\put(80,210){\oval(20,20)[r]}
\put(80,195){\oval(10,10)[l]}
\put(80,130){\oval(60,120)[tr]}
\put(100,150){\oval(40,100)[tr]}
\put(110,130){\line(1,2){10}}
\put(100,200){\line(0,1){20}}

\put(120,20){\line(0,1){40}}
\put(120,60){\line(-1,1){40}}
\put(80,100){\line(2,1){80}}
\put(160,140){\line(-1,2){40}}

\put(140,200){\circle{10}}

\put(130,100){\framebox(10,10){}}

\put(170,20){\oval(60,160)[tl]}
\put(170,20){\oval(20,120)[tl]}
\put(170,80){\oval(20,40)[tr]}
\put(175,80){\oval(10,10)[b]}

\put(180,20){\line(-1,1){10}}
\put(170,30){\line(1,0){10}}
\put(180,30){\line(-1,2){10}}
\put(180,40){\line(-1,1){10}}
\put(180,40){\line(0,1){20}}
\put(190,40){\line(-1,2){10}}
\put(190,40){\line(1,1){10}}
\put(190,30){\line(1,2){10}}
\put(190,30){\line(1,0){20}}
\put(200,20){\line(1,1){10}}

\put(220,20){\line(-2,5){80}}

\put(170,220){\oval(20,40)[bl]}
\put(170,190){\oval(20,20)[tr]}
\put(190,190){\oval(20,20)[b]}
\put(210,190){\oval(20,20)[tl]}
\put(210,220){\oval(20,40)[br]}

\put(190,200){\circle{10}}

\put(190,220){\oval(20,20)[b]}

\end{picture}
\end{center}

We say that two {\it n}-diagrams $D_{1}$ in $R_{n,a}$
and $D_{2}$ in $R_{n,b}$ are {\it equivalent}, and write
$D_1 \simeq D_2$, iff there is a homeomorphism
$h:D_1 \rightarrow D_2$ such that $h(i,0)=(i,0)$ and $h(i,a)=(i,b)$.
It is clear that this defines indeed an equivalence relation between
{\it n}-diagrams.

Equivalence classes of {\it n}-diagrams are the sort of
structure that make what Kauffman in \cite[pp. 440-441]{K4}
calls a {\it connection monoid}. Kauffman's
{\it diagram monoids} of \cite[p. 440]{K4} are generated from
equivalence classes of particular {\it n}-diagrams, such as
we consider in Section 3.
It simplifies matters that equivalence of {\it n}-diagrams
can be defined in terms
of homeomorphisms, as it was done here, rather than in terms of
ambient isotopies, as with other sorts of tangles, which involve crossings.

If $i$ stands for $(i,a)$ and $-i$ stands for $(i,0)$, we may identify
the end points of each thread of an {\it n}-diagram in $R_{n,a}$ by a
pair of integers in $[-n,n]- \{ 0 \}$. Then it is easy to see the following.

\begin{r}
The {\it n}-diagrams $D_1$ and $D_2$ are equivalent iff   
\vspace{0.1cm}

$(i)$ the end points of the threads in $D_1$ are identified
with the same pairs of integers as the end points of the threads in
$D_2$, and

$(ii)$ $D_1$ and $D_2$ have the same number of circular components.

\end{r}

Let us say that the closed interval
$[x,y]= \{ r\in {\bf R} \;|\; x \le r \le y \}$
is {\it proper} iff $x < y$. We say that $[a,b]$ {\it encloses} $[c,d]$
iff
$a<c$ and $d<b$. We may then identify the equivalence class of
an {\it n}-diagram
by a set of $n$ proper intervals in  $[-n,n]- \{ 0 \}$ with end points in
{\bf Z} such that any two distinct intervals in the set are either disjoint
or one of these two intervals encloses the other; we must moreover specify a
natural number that stands for the number of circular components in the
{\it n}-diagrams of our equivalence class. Each interval with
end points in {\bf Z}
may, of course, be identified with a pair of integers. So we may identify
the equivalence class of an {\it n}-diagram by the pair $(\Theta , l)$
where $\Theta$ is a set of $n$ pairs of
integers in $[-n,n]- \{ 0 \}$ that
satisfies the conditions stated above in terms of intervals and $l$ is the
number of circular components.

{\it Parenthetical words} are finite sequences of the symbols ( and )
defined inductively as follows:

\vspace{0.1cm}

the empty word is a parenthetical word;

if $\alpha$ is a parenthetical word, then $( \alpha )$ is
a parenthetical word;

if $\alpha$ and $\beta$ are   parenthetical words, then
$\alpha \beta$ is a parenthetical word.
\vspace{0.1cm}

\noindent It is clear that parenthetical words with $2n$ symbols are in one-to-one
correspondence with equivalence classes of {\it n}-diagrams without
circular components.

Let an {\it o-monoid}
be a monoid with an arbitrary unary operation $o$, and consider the free
{\it o}-monoid $\cal M$ generated by the empty set of generators. The
free {\it o}-monoid
$\cal M$ is isomorphic to parenthetical words where the unit is the
empty word, multiplication is concatenation, and the operation $o$ is
putting in
parentheses. So equivalence classes of {\it n}-diagrams
without circular components
may be conceived as elements of $\cal M$.

The set of equivalence classes of {\it n}-diagrams is endowed with the
structure
of a monoid in the following manner. Let the
{\it unit} {\it n}-diagram $I$
be $\{ (i,y)\;|\;i \in {\bf N} \cap [1,n] \;\mbox {and}\; y \in [0,1]\}$
in $R_{n,1}$. So $I$ has no circular component, and all of its threads 
are vertical transversal threads. We draw $I$ as follows:

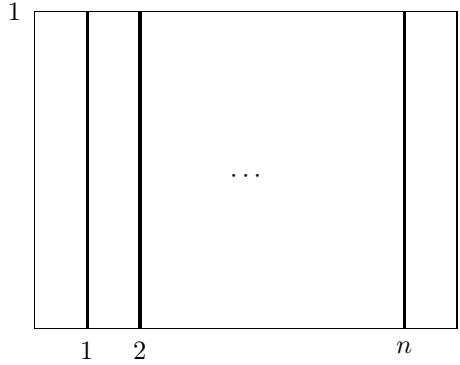
\begin{center}
\begin{picture}(160,160)
{\linethickness{0.05pt}
\put(0,20){\line(1,0){160}}
\put(0,20){\line(0,1){120}}
\put(0,140){\line(1,0){160}}
\put(160,20){\line(0,1){120}}}

\put(-5,140){\makebox(0,0)[r]{$1$}}
\put(20,15){\makebox(0,0)[t]{$1$}}
\put(40,15){\makebox(0,0)[t]{$2$}}
\put(140,15){\makebox(0,0)[t]{$n$}}

\thicklines
\put(20,20){\line(0,1){120}}
\put(40,20){\line(0,1){120}}
\put(140,20){\line(0,1){120}}

\put(80,80){\makebox(0,0)[t]{$\cdots$}}

\end{picture}
\end{center}

For two {\it n}-diagrams $D_1$ in $R_{n,a}$ and $D_2$ in $R_{n,b}$
let the {\it composition} of $D_1$ and $D_2$ be defined as follows:
$$D_2 \circ  D_1 = \{ (x, y+b) \;|\; (x,y) \in D_1 \} \cup D_2.$$

\noindent It is easy to see that $D_2 \circ D_1$ is an {\it n}-diagram in
$R_{n,a+b}$.

Let $D_i$ be an {\it n}-diagram in $R_{n,a_i}$, and suppose 
$D_1 \simeq D_3$ with the homeomorphism
$h_1:D_1 \rightarrow D_3$, and
$D_2 \simeq D_4$ with the homeomorphism
$h_2:D_2 \rightarrow D_4$.
Then $D_2 \circ D_1 \simeq D_4 \circ D_3$ with the
homeomorphism $h:D_2 \circ D_1 \rightarrow D_4 \circ D_3$ defined as
follows. For $p^1$ the first and $p^2$ the second projection, let

$$ h(x,y)= \left \{ \begin {array}{ll}
(p^1(h_1(x, y-a_2)), p^2(h_1(x, y-a_2))+a_4), & \mbox { if } y>a_2 \\
h_2(x, y), & \mbox { if } y \le a_2.
\end{array} \right. $$

So the composition $\circ$ defines an operation on equivalence
classes of {\it n}-diagrams. We can then establish that

\vspace{0.1cm}

$(1) \; \; I \circ D \simeq D \circ I \simeq D,$

\vspace{0.1cm}

$(2) \; \; D_3 \circ (D_2 \circ D_1) \simeq (D_3 \circ D_2) \circ D_1.$
\vspace{0.1cm}

The equivalences of (1) follow from the fact that 
$I \circ D$, $D \circ I$ and $D$ have the same number of
circular components, because $I$ has no circular component,
and that their threads
may be identified with the same pairs of integers,
because all the threads of $I$ are vertical transversal threads.
Then we apply Remark 3.
For the equivalence (2), it is clear that $D_3 \circ (D_2 \circ D_1)$
is actually identical to  $(D_3 \circ D_2) \circ D_1.$

Let $a(n)$ be $\mbox{\rm max}\{5,(n-1)(n-2)/2\}$.
(The precise value of $a(n)$ could be varied, and we
have chosen one of infinitely many possibilities.)
We say that an {\it n}-diagram is {\it normal} iff
it is in $R_{n,a(n)}$,
each of its transversal threads is a straight line segment, and each of
its cups and caps is a semicircle. It is clear that two
normal {\it n}-diagrams without
circular components are equivalent iff they are equal.
(To make this equivalence
hold even in the presence of circular components we can find a
definite place
in $R_{n,a(n)}$ for an arbitrary finite number of circles that
will not intersect with each other and with other components.)

Every {\it n}-diagram $D$ is equivalent to a normal
{\it n}-diagram, which is a 
handy representative for the equivalence class of $D$.
A normal 11-diagram equivalent to the 11-diagram 
from the beginning of this section looks as follows:
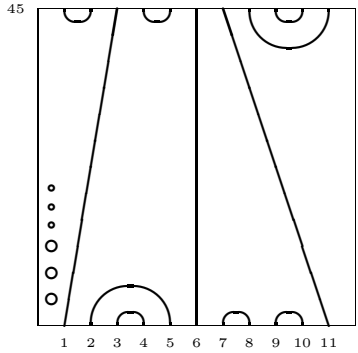
\vspace{-0.2cm}

\begin{center}
\begin{picture}(120,160)
{\linethickness{0.05pt}
\put(0,20){\line(1,0){120}}
\put(0,20){\line(0,1){120}}
\put(0,140){\line(1,0){120}}
\put(120,20){\line(0,1){120}}}

\put(-5,140){\makebox(0,0)[r]{\tiny$45$}}
\put(10,15){\makebox(0,0)[t]{\tiny$1$}}
\put(20,15){\makebox(0,0)[t]{\tiny$2$}}
\put(30,15){\makebox(0,0)[t]{\tiny$3$}}
\put(40,15){\makebox(0,0)[t]{\tiny$4$}}
\put(50,15){\makebox(0,0)[t]{\tiny$5$}}
\put(60,15){\makebox(0,0)[t]{\tiny$6$}}
\put(70,15){\makebox(0,0)[t]{\tiny$7$}}
\put(80,15){\makebox(0,0)[t]{\tiny$8$}}
\put(90,15){\makebox(0,0)[t]{\tiny$9$}}
\put(100,15){\makebox(0,0)[t]{\tiny$10$}}
\put(110,15){\makebox(0,0)[t]{\tiny$11$}}

\thicklines
\put(10,20){\line(1,6){20}}
\put(60,20){\line(0,1){120}}
\put(110,20){\line(-1,3){40}}

\put(5,30){\circle{4}}
\put(5,40){\circle{4}}
\put(5,50){\circle{4}}
\put(5,58){\circle{2}}
\put(5,65){\circle{2}}
\put(5,72){\circle{2}}

\put(35,20){\oval(30,30)[t]}
\put(35,20){\oval(10,10)[t]}
\put(75,20){\oval(10,10)[t]}
\put(95,20){\oval(10,10)[t]}

\put(15,140){\oval(10,10)[b]}
\put(45,140){\oval(10,10)[b]}
\put(95,140){\oval(10,10)[b]}
\put(95,140){\oval(30,30)[b]}

\end{picture}
\end{center}

\vspace{-0.2cm}

\section{Generating {\it n}-diagrams}
\noindent For $i \in {\bf N} \cap [1,n-1]$, the {\it diapsidal}
{\it n}-diagram $H^i$
is the normal {\it n}-diagram without circular components, with a single
cup with
the end points $(i,a(n))$ and $(i+1,a(n))$, and a single cap
with the end points $(i,0)$ and $(i+1,0)$; all the other threads are
transversal
threads orthogonal to the $x$ axis. A diapsidal {\it n}-diagram $H^i$
looks as follows:
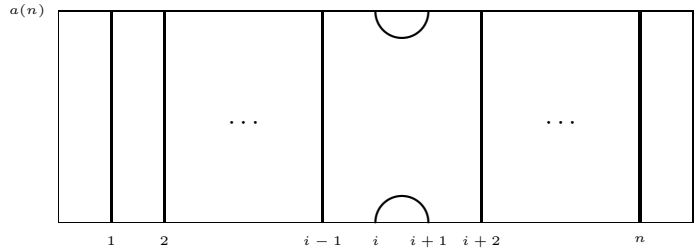
\vspace{-0.2cm}

\begin{center}
\begin{picture}(240,120)
{\linethickness{0.05pt}
\put(0,20){\line(1,0){240}}
\put(0,20){\line(0,1){80}}
\put(0,100){\line(1,0){240}}
\put(240,20){\line(0,1){80}}}

\put(-5,100){\makebox(0,0)[r]{\tiny$a(n)$}}
\put(20,15){\makebox(0,0)[t]{\tiny$1$}}
\put(40,15){\makebox(0,0)[t]{\tiny$2$}}
\put(100,15){\makebox(0,0)[t]{\tiny$i-1$}}
\put(120,15){\makebox(0,0)[t]{\tiny$i$}}
\put(140,15){\makebox(0,0)[t]{\tiny$i+1$}}
\put(160,15){\makebox(0,0)[t]{\tiny$i+2$}}
\put(220,15){\makebox(0,0)[t]{\tiny$n$}}

\thicklines
\put(20,20){\line(0,1){80}}
\put(40,20){\line(0,1){80}}
\put(100,20){\line(0,1){80}}
\put(160,20){\line(0,1){80}}
\put(220,20){\line(0,1){80}}

\put(130,20){\oval(20,20)[t]}
\put(130,100){\oval(20,20)[b]}

\put(70,60){\makebox(0,0)[t]{$\cdots$}}
\put(190,60){\makebox(0,0)[t]{$\cdots$}}

\end{picture}
\end{center}

The {\it circular} {\it n}-diagram $C$ is the {\it n}-diagram that
differs from the unit
{\it n}-diagram $I$ by having a single circular component,
which for the sake of definiteness
we may choose to be a circle of radius, let us say, 1/4,
with centre $(1/2,1/2)$.

If the end points of a thread of an {\it n}-diagram are $(i,x)$
and $(j,y)$, where
$x,y \in \{0,a \}$, let us say that this thread {\it covers} a pair of
natural
numbers $(m,l)$, where $1 \le m<l \le n$, iff $min \{ i,j \} \le m$ and
$l \le max \{ i,j \}$.
Then we can establish the following.

\begin{r}
In every {\it n}-diagram, every pair $(m,m+1)$, where $1 \le m<n$,
is covered by an even number of threads.
\end{r}

\noindent {\it Proof.} For an {\it n}-diagram $D$ the cardinality of the set
$P= \{ (i,x) \in t(D) \cup b(D)\;|\; i \le m \}$ is $2m$.
Every thread of $D$ that covers $(m,m+1)$ has a single end point in $P$,
and other threads of $D$ have 0 or 2 end points in $P$. Since $P$ is even,
the remark follows. {\it q.e.d.}
\vspace{0.2cm}

If the end points of a thread of an {\it n}-diagram are $(i,x)$
and $(j,y)$, let
the {\it span of the thread} be $|i-j|$. Let the {\it span} $\sigma (D)$
{\it of an n-diagram} $D$ be the sum of the spans of all the
threads of $D$. Remark 4 entails that the span of an
{\it n}-diagram is an even number greater than or equal to 0; the
span of $I$ is 0. It is clear that equivalent {\it n}-diagrams have the
same span. It is also easy to see that the following holds.

\begin{r}
If an {\it n}-diagram has cups, then it must have at least one cup
whose span is $1$. The same holds for caps.
\end{r}

We can now prove the following lemma.

\begin{lm}
Every {\it n}-diagram is equivalent to an {\it n}-diagram generated
from $I$, $C$ and the diapsidal {\it n}-diagrams
$H^i$, for $1 \le i \le n-1$, with the operation of composition $\circ$.
\end{lm}
{\it Proof.} Take an arbitrary {\it n}-diagram $D$ in $R_{n,a}$, and
let $D$ have $l \ge 1$ circular components.
Let $C^1$ be $C$, and let $C^{l+1}$
be $C^l \circ C$. It is clear that $D$ is equivalent to an {\it n}-diagram
$C^l \circ D_1$ in  $R_{n,a+l}$ with $D_1$ an
{\it n}-diagram in $R_{n,a}$ without circular components.
If $l=0$, then $D_1$ is $D$. It simplifies matters if we assume that
$D_1$ is a normal {\it n}-diagram, but this is not essential. 

Then we proceed by induction on  $\sigma (D_1)$. If $\sigma (D_1)=0$,
then $D_1 \simeq I$, and the lemma holds. Suppose  
$\sigma (D_1)>0$. Then there must be a cup or a cap in $D_1$,
and by Remarks 2 and 5, there must be at least one cup of $D_1$
whose span is 1.
The end points of such cups are of the form $(i,a)$ and $(i+1,a)$.
Then select among
all these cups that one where $i$ is the greatest number; let that number
be $j$, and let us call the cup we have selected $\upsilon_j$.
By Remark 4,
there must be at least one other thread of $D_1$,
different from $\upsilon_j$,
that also covers $(j,j+1)$. We have to consider four cases, which exclude
each other:

\begin{itemize}
\item[$(1)$] {\rm $(j,j+1)$ is covered by a cup $\xi$ of $D_1$ different
from $\upsilon_j$, whose end points are $(p,a)$ and $(q,a)$ with $p<q$;}
\item[$(2.1)$] $(j,j+1)$ {\rm is not covered by a cup of $D_1$ different
from $\upsilon_j$, but it is covered by a transversal thread
$\xi$ of $D_1$ whose end points are $(p,a)$ and $(q,0)$ with $p<q$;}
\item[$(2.2)$] {\rm same as case (2.1) save that the end points of
$\xi$ are $(p,0)$ and $(q,a)$ with $p<q$;}
\item[$(3)$] $(j,j+1)$ {\rm is covered neither by a cup of $D_1$ different
from $\upsilon_j$, nor by a transversal thread
of $D_1$, but it is covered by a cap $\xi$ of $D_1$, whose end points
are $(p,0)$ and $(q,0)$ with $p<q$.}
\end{itemize}

In cases (1) and (2.1) we select among the threads $\xi$ mentioned
the one where $p$ is maximal, and in cases
(2.2) and (3) we select the $\xi$ where $p$ is minimal.
(We obtain the same result if in cases (1) and (2.2) we take $q$ minimal,
while in (2.1) and (3) we take $q$ maximal.)

We build out of the {\it n}-diagram $D_1$ a new {\it n}-diagram $D_2$ in
$R_{n,a}$ by replacing the thread $\upsilon_j$ and the selected thread
$\xi$, whose end points are $(p,x)$ and $(q,y)$, where $x,y \in \{ 0,a \}$,
with two new threads: one whose end points are $(p,x)$ and $(j,a)$,
and the other
whose end points are $(j+1,a)$ and $(q,y)$. We can easily check that
$D_1 \simeq D_2 \circ H^j$. This is clear from the following picture:

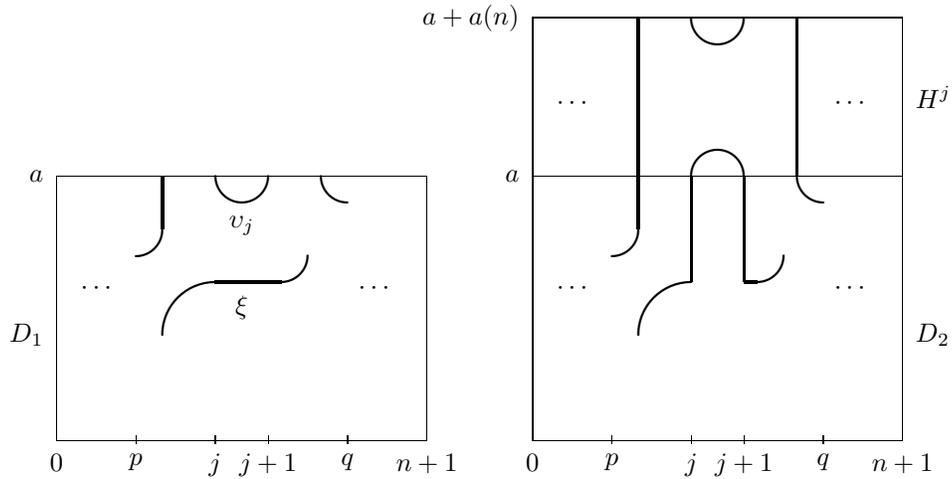
\begin{center}
\begin{picture}(320,200)
{\linethickness{0.05pt}
\put(0,20){\line(1,0){140}}
\put(0,20){\line(0,1){100}}
\put(0,120){\line(1,0){140}}
\put(140,20){\line(0,1){100}}

\put(180,20){\line(1,0){140}}
\put(180,20){\line(0,1){160}}
\put(180,180){\line(1,0){140}}
\put(320,20){\line(0,1){160}}
\put(30,19){\line(0,1){3}}
\put(60,19){\line(0,1){3}}
\put(80,19){\line(0,1){3}}
\put(110,19){\line(0,1){3}}
\put(210,19){\line(0,1){3}}
\put(240,19){\line(0,1){3}}
\put(260,19){\line(0,1){3}}
\put(290,19){\line(0,1){3}}
}

\put(-5,120){\makebox(0,0)[r]{$a$}}
\put(0,15){\makebox(0,0)[t]{$0$}}
\put(30,15){\makebox(0,0)[t]{$p$}}
\put(60,15){\makebox(0,0)[t]{$j$}}
\put(80,15){\makebox(0,0)[t]{$j+1$}}
\put(110,15){\makebox(0,0)[t]{$q$}}
\put(140,15){\makebox(0,0)[t]{$n+1$}}

\put(180,15){\makebox(0,0)[t]{$0$}}
\put(175,120){\makebox(0,0)[r]{$a$}}
\put(175,180){\makebox(0,0)[r]{$a+a(n)$}}
\put(210,15){\makebox(0,0)[t]{$p$}}
\put(240,15){\makebox(0,0)[t]{$j$}}
\put(260,15){\makebox(0,0)[t]{$j+1$}}
\put(290,15){\makebox(0,0)[t]{$q$}}
\put(320,15){\makebox(0,0)[t]{$n+1$}}
\put(70,105){\makebox(0,0)[t]{$\upsilon_j$}}
\put(70,75){\makebox(0,0)[t]{$\xi$}}

\thicklines
\put(30,120){\oval(20,60)[br]}
\put(70,120){\oval(20,20)[b]}
\put(110,120){\oval(20,20)[bl]}
\put(80,60){\oval(80,40)[tl]}
\put(80,90){\oval(30,20)[br]}

\put(15,80){\makebox(0,0)[t]{$\cdots$}}
\put(120,80){\makebox(0,0)[t]{$\cdots$}}

\put(-5,60){\makebox(0,0)[r]{$D_1$}}

\put(210,120){\oval(20,60)[br]}
\put(290,120){\oval(20,20)[bl]}
\put(240,60){\oval(40,40)[tl]}
\put(260,90){\oval(30,20)[br]}
\put(240,80){\line(0,1){40}}
\put(260,80){\line(0,1){40}}

\put(220,120){\line(0,1){60}}
\put(280,120){\line(0,1){60}}

\put(250,120){\oval(20,20)[t]}
\put(250,180){\oval(20,20)[b]}

\put(195,80){\makebox(0,0)[t]{$\cdots$}}
\put(300,80){\makebox(0,0)[t]{$\cdots$}}
\put(195,150){\makebox(0,0)[t]{$\cdots$}}
\put(300,150){\makebox(0,0)[t]{$\cdots$}}

\put(325,60){\makebox(0,0)[l]{$D_2$}}
\put(325,150){\makebox(0,0)[l]{$H^j$}}

{\linethickness{0.05pt}
\put(180,120){\line(1,0){140}}}

\end{picture}
\end{center}

Neither of the new threads of $D_2$ that have replaced $\upsilon_j$
and $\xi$ covers $(j,j+1)$, and $\sigma (D_1)= \sigma (D_2)+2$.
So, by the induction hypothesis, $D_2$ is equivalent to an
{\it n}-diagram $D_3$ generated from $I$, $C$ and $H^i$ with $\circ$,
and since
$D_1 \simeq D_3 \circ H^j$, this proves the lemma. {\it q.e.d.}
\vspace{0.2cm}

Note that we need not require in this proof that in $\upsilon_j$
the number $j$ should be the greatest number $i$ for cups with end points
$(i,a)$ and $(i+1,a)$. The proof would go through without making this choice.
But with this choice we shall end up with a composition of {\it n}-diagrams
that
corresponds exactly to a term of ${\cal K}_{n}$ in Jones normal form.

\section{${\cal K}_{n}$ is the monoid of {\it n}-diagrams}

\noindent Let ${\cal D}_{n}$ be the set of {\it n}-diagrams. We define a map 
$\delta :{\cal K}_{n} \rightarrow {\cal D}_{n}$ as follows:

\[\begin{array}{lcl}
\delta (h^i) & = & H^{i},\\[0.1cm]
\delta (c) & = & C,\\[0.1cm]
\delta ({\bf 1}) & = & I,\\[0.1cm]
\delta (tu) & = & \delta (t) \circ \delta (u).
\end{array}\]

We can then prove the following.

\begin{lm}
If $t=u$ in  ${\cal K}_{n}$, then $\delta (t) \simeq \delta (u)$. 
\end{lm}
{\it Proof.} We have already verified in Section 2 that we have
replacement of equivalents, and that the equations (1) and (2)
of the axiomatization of ${\cal K}_{n}$ are satisfied for $I$ and
$\circ$. It just remains to verify $(h1), (h2), (hc1)$ and $(hc2)$,
which is quite straightforward. {\it q.e.d.}
\vspace{0.2cm}

Let [${\cal D}_n$] be the set of equivalence classes $[D]$ of the
{\it n}-diagrams $D$.
This set is a monoid whose unit is $[I]$ and whose multiplication is
defined by taking that $[D_1] [D_2]$ is $[D_1 \circ D_2]$.
If $\delta ':{\cal K}_{n} \rightarrow [{\cal D}_{n}]$ is defined as $\delta$
on the generators of ${\cal K}_{n}$, save that $H^i$ and $C$ are replaced
by $[H^i]$ and $[C]$, and if $\delta '(\bf 1)$ is $[I]$ and 
$\delta '(tu)= \delta '(t) \delta '(u)$, then Lemma 3 guarantees that $\delta '$ is a
homomorphism. Lemma 2 guarantees that $\delta '$ is onto. To establish that
$\delta '$ is an isomorphism it remains only to show that $\delta '$ is one-one.

Let a transversal thread in an {\it n}-diagram be called {\it falling} iff
its end points are $(i,a)$ and $(j,0)$ with $i<j$. If the end points of
a thread of an {\it n}-diagram $D$ are
$(i,a)$ and $(j,x)$ with $x \in \{ 0,a \}$ and $i<j$, then we say that
$(i,a)$ is a {\it top slope point} of $D$. If the end points of a
thread of $D$ are $(i,x)$ and $(j,0)$ with $x \in \{ 0,a \}$ and $i<j$,
then we say that $(j,0)$ is a {\it bottom slope point} of $D$. Each cup
has a single top slope point, each cap has a single bottom slope point,
and each falling transversal thread has one top slope point and one
bottom slope point. Other transversal threads have no slope points.
So we can ascertain the following.

\begin{r}
In every {\it n}-diagram the number of top slope points is equal
to the number  of bottom slope points.
\end{r}

\noindent Remember that by Remark 2 the number of cups is equal to the
number of caps.

Let $(a_1,a), \ldots , (a_k,a)$ be the sequence of all top slope points
of an  {\it n}-diagram $D$, ordered so that $a_1< \ldots <a_k$, and let
$(b_1+1,0), \ldots ,(b_k+1,0)$ be the sequence of all bottom slope
points of $D$,
ordered so that $b_1< \ldots <b_k$ (as we just saw with Remark 6,
these sequences must be of equal length). Then let $T_D$ be the sequence
of natural numbers $a_1, \ldots , a_k$ and $B_D$ the sequence of natural
numbers  $b_1, \ldots , b_k$.

\begin{r}
The sequence $T_{\delta (h^{[i,j]})}$ has a single member $j$  
and the sequence $B_{\delta (h^{[i,j]})}$ has a single member $i$.
\end{r}

\noindent This is clear from the {\it n}-diagram $\delta (h^{[i,j]})$,
which is equivalent to an {\it n}-diagram of the following form:

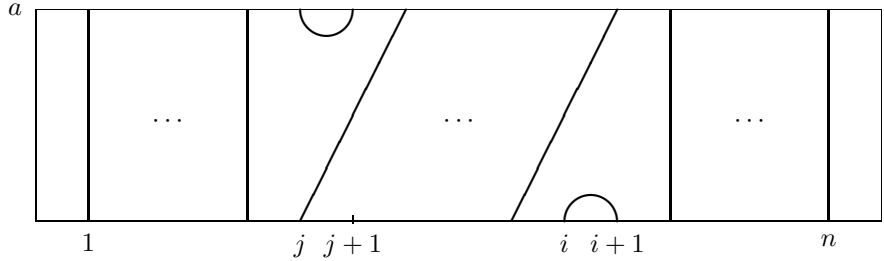
\begin{center}
\begin{picture}(320,120)
{\linethickness{0.05pt}
\put(0,20){\line(1,0){320}}
\put(0,20){\line(0,1){80}}
\put(0,100){\line(1,0){320}}
\put(320,20){\line(0,1){80}}}
\put(120,19){\line(0,1){3}}

\put(-5,100){\makebox(0,0)[r]{$a$}}
\put(20,15){\makebox(0,0)[t]{$1$}}
\put(100,15){\makebox(0,0)[t]{$j$}}
\put(120,15){\makebox(0,0)[t]{$j+1$}}
\put(200,15){\makebox(0,0)[t]{$i$}}
\put(220,15){\makebox(0,0)[t]{$i+1$}}
\put(300,15){\makebox(0,0)[t]{$n$}}

\thicklines
\put(20,20){\line(0,1){80}}
\put(80,20){\line(0,1){80}}
\put(240,20){\line(0,1){80}}
\put(300,20){\line(0,1){80}}

\put(100,20){\line(1,2){40}}
\put(180,20){\line(1,2){40}}

\put(110,100){\oval(20,20)[b]}
\put(210,20){\oval(20,20)[t]}

\put(50,60){\makebox(0,0)[t]{$\cdots$}}
\put(160,60){\makebox(0,0)[t]{$\cdots$}}
\put(270,60){\makebox(0,0)[t]{$\cdots$}}

\end{picture}
\end{center}

\noindent provided $0<j<i<n$ (in other cases we simplify this
diagram by omitting some
transversal threads).

\begin{r}
If $D_1 \simeq D_2$, then $T_{D_1} = T_{D_2}$ and
$B_{D_1} = B_{D_2}.$
\end{r}

\noindent This follows from Remark 3.

Then we can prove the following lemmata.

\begin{lm}
If $t$ is the term
$h^{[b_1,a_1]} \ldots h^{[b_k,a_k]}$ with
$a_1< \ldots <a_k$ and $b_1< \ldots <b_k$,
then $T_{\delta (t)}$ is $a_1, \ldots , a_k$ and
$B_{\delta (t)}$ is $b_1, \ldots , b_k$.
\end{lm}
{\it Proof.} We proceed by induction on $k$. If $k=1$,
we use Remark 7. If $k>1$, then, by the induction hypothesis, the lemma
has been established for the term
$h^{[b_1,a_1]} \ldots h^{[b_{k-1},a_{k-1}]}$,
which we call $t'$. So $T_{\delta (t')}$ is $a_1, \ldots , a_{k-1}$.

Since in $\delta (h^{[b_k,a_k]})$ every point in the top with the first
coordinate $i<a_k$ is the end point of a vertical transversal thread,
and since $\delta (t)= \delta (t') \circ \delta (h^{[b_k,a_k]})$, the beginning
of the sequence $T_{\delta (t)}$ must be $a_1, \ldots , a_{k-1}$. To this
sequence we have to add $a_k$ because $\delta (t)$ inherits the cup of
$\delta (h^{[b_k,a_k]})$. This shows immediately that $a_k+1$ is not in 
$T_{\delta (t)}$. It remains to show that for no $i \ge a_k+2$ we can have in
$\delta (t)$ a top slope point with the first coordinate $i$.

If $i>b_k+1$, then every point in the top with the first coordinate $i$ is
the end point of a vertical transversal thread in both $\delta (h^{[b_k,a_k]})$
and $\delta (t')$. So $i$ is not in $T_{\delta (t)}$. It remains to consider
$i$ for $a_k+2 \le i \le b_k+1$. Every point in the top with
this first coordinate $i$ is the end point of a transversal thread in
$\delta (h^{[b_k,a_k]})$ whose other end point is $(i-2,0)$. If $i$ were to
be added to  $T_{\delta (t)}$, the number $i-2$ would be in $T_{\delta (t')}$,
but this contradicts the fact that $T_{\delta (t')}$ ends with $a_{k-1}$. So
$T_{\delta (t)}$ is $a_1, \ldots , a_k$.

To show that  $B_{\delta (t)}$ is $b_1, \ldots , b_k$ we reason analogously by
applying the induction hypothesis to $h^{[b_2,a_2]} \ldots h^{[b_k,a_k]}$.
{\it q.e.d.}
\vspace{0.2cm}

\begin{lm}
If $t$ and $u$ are terms of ${\cal K}_{n}$ in Jones normal form and
 $\delta (t) \simeq \delta (u)$, then $t$ and $u$ are the same term.
\end{lm}
{\it Proof.} Let $t$ be $c^l h^{[b_1,a_1]} \ldots h^{[b_k,a_k]}$ and
let $u$ be $c^j h^{[d_1, c_1]} \ldots h^{[d_m, c_m]}$. If $l \neq j$, then 
$\delta (t)$ is not equivalent to $\delta (u)$ by Remark 3, because
$\delta (t)$ and $\delta (u)$ have different numbers of circular components.
If $a_1, \ldots , a_k$ is different from  $c_1, \ldots , c_m$, or
$b_1, \ldots , b_k$
is different from $d_1, \ldots , d_m$, then $\delta (t)$ is not equivalent to
$\delta (u)$ by Lemma 4 and Remark 8. {\it q.e.d.}
\vspace{0.2cm}

\begin{lm}
If $\delta (t) \simeq \delta (u)$, then $t=u$ in ${\cal K}_{n}$.
\end{lm}
{\it Proof.} Suppose $t=u$ doesn't hold in ${\cal K}_{n}$. Let
$t'$ and $u'$ be terms in the Jones normal form such that
$t=t'$ and $u=u'$ in ${\cal K}_{n}$. Such terms exist according
to Lemma 1.
Then $t'$ and $u'$ must be different  
terms; otherwise $t=u$ would hold in ${\cal K}_{n}$. By Lemma 5 we have that
$\delta (t')$ is not equivalent to $\delta (u')$. Since by Lemma 3
we have  $\delta (t) \simeq \delta (t')$ and $\delta (u) \simeq \delta (u')$, we can
 conclude  that $\delta (t)$ is not equivalent to $\delta (u)$. {\it q.e.d.}
\vspace{0.2cm}

Lemma 6 guarantees that $\delta ':{\cal K}_{n} \rightarrow [{\cal D}_{n}]$
is one-one, and so $\delta '$ is an isomorphism.

We can now conclude that for every term $t$ of ${\cal K}_{n}$ there is a
{\it unique} term $t'$ in Jones normal form such that $t=t'$ in
${\cal K}_{n}$. Otherwise, if $t$ were equal in ${\cal K}_{n}$ to two
different terms $t'$ and $t''$ in Jones normal form, by Lemma 3
we would have
$\delta (t) \simeq \delta (t')$ and $\delta (t) \simeq \delta (t'')$, and hence
also $\delta (t') \simeq \delta (t'')$. But this contradicts Lemma 5.

This solves the word problem for Kauffman monoids. To check whether $t=u$ in
${\cal K}_{n}$ just reduce $t$ and $u$ to Jones normal form, according to the procedure
of the proof of Lemma 1, and then check whether the normal forms obtained are
equal. However, to reduce a term of ${\cal K}_{n}$ to Jones normal
form, now that we have established that $\delta '$ is an isomorphism, we can
proceed more efficiently with {\it n}-diagrams than with the syntactical
method of the proof of Lemma 1.

From Remark 1 in Section 1 it follows that in
the {\it n}-diagram $\delta (t)$
of a term $t$ in Jones normal form, or in any alternative
reduced normal form,
such as those envisaged after Remark 1, we will have no threads with
bulges like the following:

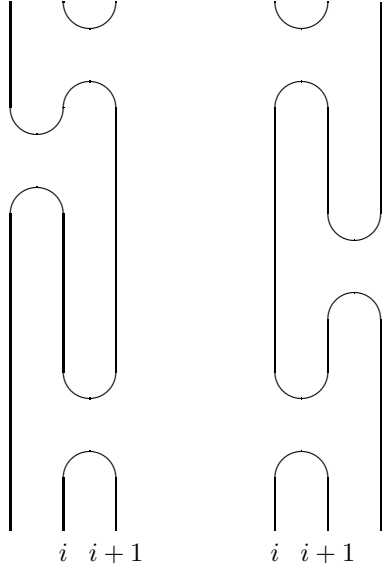
\begin{center}
\begin{picture}(140,230)

\put(20,15){\makebox(0,0)[t]{$i$}}
\put(40,15){\makebox(0,0)[t]{$i+1$}}

\put(100,15){\makebox(0,0)[t]{$i$}}
\put(120,15){\makebox(0,0)[t]{$i+1$}}

\put(0,20){\line(0,1){120}}
\put(0,180){\line(0,1){40}}
\put(20,20){\line(0,1){20}}
\put(40,20){\line(0,1){20}}
\put(20,80){\line(0,1){60}}
\put(40,80){\line(0,1){100}}

\put(30,40){\oval(20,20)[t]}
\put(30,80){\oval(20,20)[b]}

\put(30,180){\oval(20,20)[t]}
\put(30,220){\oval(20,20)[b]}

\put(10,140){\oval(20,20)[t]}
\put(10,180){\oval(20,20)[b]}

\put(100,20){\line(0,1){20}}
\put(100,80){\line(0,1){100}}
\put(120,20){\line(0,1){20}}
\put(140,20){\line(0,1){80}}
\put(120,80){\line(0,1){20}}
\put(140,140){\line(0,1){80}}
\put(120,140){\line(0,1){40}}

\put(130,100){\oval(20,20)[t]}
\put(130,140){\oval(20,20)[b]}

\put(110,40){\oval(20,20)[t]}
\put(110,80){\oval(20,20)[b]}

\put(110,180){\oval(20,20)[t]}
\put(110,220){\oval(20,20)[b]}

\end{picture}
\end{center}

\noindent Reduction to normal form involves getting rid of these bulges,
and this is done more easily diagrammatically than syntactically.
To reduce a term $t$ to Jones normal form we first draw the {\it n}-diagram
$\delta (t)$. Then we replace $\delta (t)$ by a normal {\it n}-diagram $D$
such that $D \simeq \delta (t)$. This is where bulges get eliminated.

Our proof of Lemma 2 in Section 3 then gives a procedure for
building out of $D$ the term $t'$ in Jones normal form
such that $\delta (t') \simeq D$, Lemma 6
guaranteeing that $t'$ is equal in ${\cal K}_{n}$ to the original $t$.
We could apply this procedure to $\delta (t)$ directly, but it is easier to
apply it to the normal {\it n}-diagram $D$.

There is, however, a procedure handier than that,
which also yields a term in Jones normal form out of a normal
{\it n}-diagram.
This procedure (suggested by Figure 16 of
\cite[p. 434]{K4}, and detailed in \cite[Section 6]{K8}, in mirror
image)
is illustrated in the following picture, based on
the example from the end of Section 2:

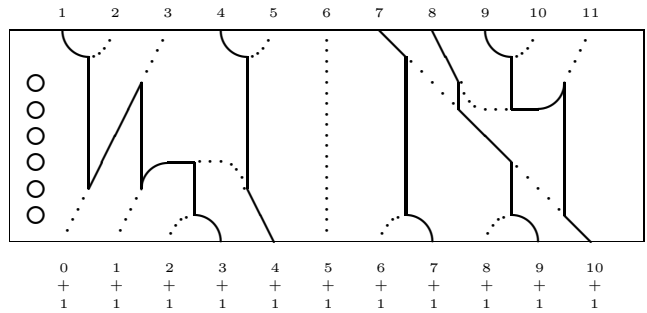
\begin{center}
\begin{picture}(240,120)
{\linethickness{0.05pt}
\put(0,20){\line(1,0){240}}
\put(0,20){\line(0,1){80}}
\put(0,100){\line(1,0){240}}
\put(240,20){\line(0,1){80}}}

\put(18,-5){\shortstack{\tiny$0$\\ \tiny$+$\\ \tiny$1$}}
\put(38,-5){\shortstack{\tiny$1$\\ \tiny$+$\\ \tiny$1$}}
\put(58,-5){\shortstack{\tiny$2$\\ \tiny$+$\\ \tiny$1$}}
\put(78,-5){\shortstack{\tiny$3$\\ \tiny$+$\\ \tiny$1$}}
\put(98,-5){\shortstack{\tiny$4$\\ \tiny$+$\\ \tiny$1$}}
\put(118,-5){\shortstack{\tiny$5$\\ \tiny$+$\\ \tiny$1$}}
\put(138,-5){\shortstack{\tiny$6$\\ \tiny$+$\\ \tiny$1$}}
\put(158,-5){\shortstack{\tiny$7$\\ \tiny$+$\\ \tiny$1$}}
\put(178,-5){\shortstack{\tiny$8$\\ \tiny$+$\\ \tiny$1$}}
\put(198,-5){\shortstack{\tiny$9$\\ \tiny$+$\\ \tiny$1$}}
\put(218,-5){\shortstack{\tiny$10$\\ \tiny$+$\\ \tiny$1$}}

\put(20,105){\makebox(0,0)[b]{\tiny$1$}}
\put(40,105){\makebox(0,0)[b]{\tiny$2$}}
\put(60,105){\makebox(0,0)[b]{\tiny$3$}}
\put(80,105){\makebox(0,0)[b]{\tiny$4$}}
\put(100,105){\makebox(0,0)[b]{\tiny$5$}}
\put(120,105){\makebox(0,0)[b]{\tiny$6$}}
\put(140,105){\makebox(0,0)[b]{\tiny$7$}}
\put(160,105){\makebox(0,0)[b]{\tiny$8$}}
\put(180,105){\makebox(0,0)[b]{\tiny$9$}}
\put(200,105){\makebox(0,0)[b]{\tiny$10$}}
\put(220,105){\makebox(0,0)[b]{\tiny$11$}}

\thicklines

\put(10,30){\circle{6}}
\put(10,40){\circle{6}}
\put(10,50){\circle{6}}
\put(10,60){\circle{6}}
\put(10,70){\circle{6}}
\put(10,80){\circle{6}}

\put(30,100){\oval(20,20)[bl]}
\put(30,90){\line(0,-1){50}}
\put(30,40){\line(1,2){20}}
\put(50,80){\line(0,-1){40}}
\put(70,40){\oval(40,20)[tl]}
\put(70,50){\line(0,-1){20}}
\put(70,20){\oval(20,20)[tr]}

\put(22,24){\makebox(0,0){$\cdot$}}
\put(24,28){\makebox(0,0){$\cdot$}}
\put(26,32){\makebox(0,0){$\cdot$}}
\put(28,36){\makebox(0,0){$\cdot$}}

\put(52,84){\makebox(0,0){$\cdot$}}
\put(54,88){\makebox(0,0){$\cdot$}}
\put(56,92){\makebox(0,0){$\cdot$}}
\put(58,96){\makebox(0,0){$\cdot$}}

\put(32,90){\makebox(0,0){$\cdot$}}
\put(34,91){\makebox(0,0){$\cdot$}}
\put(36,93){\makebox(0,0){$\cdot$}}
\put(38,96){\makebox(0,0){$\cdot$}}

\put(42,24){\makebox(0,0){$\cdot$}}
\put(44,28){\makebox(0,0){$\cdot$}}
\put(46,32){\makebox(0,0){$\cdot$}}
\put(48,36){\makebox(0,0){$\cdot$}}

\put(61,23){\makebox(0,0){$\cdot$}}
\put(63,26){\makebox(0,0){$\cdot$}}
\put(65,28){\makebox(0,0){$\cdot$}}
\put(67,29){\makebox(0,0){$\cdot$}}

\put(72,50){\makebox(0,0){$\cdot$}}
\put(76,50){\makebox(0,0){$\cdot$}}
\put(80,50){\makebox(0,0){$\cdot$}}
\put(84,49){\makebox(0,0){$\cdot$}}
\put(86,47){\makebox(0,0){$\cdot$}}
\put(88,44){\makebox(0,0){$\cdot$}}
\put(89,40){\makebox(0,0){$\cdot$}}

\put(90,100){\oval(20,20)[bl]}

\put(90,40){\line(0,1){50}}

\put(90,40){\line(1,-2){10}}

\put(92,90){\makebox(0,0){$\cdot$}}
\put(94,91){\makebox(0,0){$\cdot$}}
\put(96,93){\makebox(0,0){$\cdot$}}
\put(98,96){\makebox(0,0){$\cdot$}}

\put(120,24){\makebox(0,0){$\cdot$}}
\put(120,28){\makebox(0,0){$\cdot$}}
\put(120,32){\makebox(0,0){$\cdot$}}
\put(120,36){\makebox(0,0){$\cdot$}}
\put(120,40){\makebox(0,0){$\cdot$}}
\put(120,44){\makebox(0,0){$\cdot$}}
\put(120,48){\makebox(0,0){$\cdot$}}
\put(120,52){\makebox(0,0){$\cdot$}}
\put(120,56){\makebox(0,0){$\cdot$}}
\put(120,60){\makebox(0,0){$\cdot$}}
\put(120,64){\makebox(0,0){$\cdot$}}
\put(120,68){\makebox(0,0){$\cdot$}}
\put(120,72){\makebox(0,0){$\cdot$}}
\put(120,76){\makebox(0,0){$\cdot$}}
\put(120,80){\makebox(0,0){$\cdot$}}
\put(120,84){\makebox(0,0){$\cdot$}}
\put(120,88){\makebox(0,0){$\cdot$}}
\put(120,92){\makebox(0,0){$\cdot$}}
\put(120,96){\makebox(0,0){$\cdot$}}

\put(140,100){\line(1,-1){10}}
\put(150,90){\line(0,-1){60}}
\put(150,20){\oval(20,20)[tr]}

\put(141,23){\makebox(0,0){$\cdot$}}
\put(143,26){\makebox(0,0){$\cdot$}}
\put(145,28){\makebox(0,0){$\cdot$}}
\put(147,29){\makebox(0,0){$\cdot$}}

\put(152,88){\makebox(0,0){$\cdot$}}
\put(156,84){\makebox(0,0){$\cdot$}}
\put(160,80){\makebox(0,0){$\cdot$}}
\put(164,76){\makebox(0,0){$\cdot$}}
\put(168,72){\makebox(0,0){$\cdot$}}

\put(160,100){\line(1,-2){10}}
\put(170,80){\line(0,-1){10}}
\put(170,70){\line(1,-1){20}}
\put(190,50){\line(0,-1){20}}
\put(190,20){\oval(20,20)[tr]}

\put(181,23){\makebox(0,0){$\cdot$}}
\put(183,26){\makebox(0,0){$\cdot$}}
\put(185,28){\makebox(0,0){$\cdot$}}
\put(187,29){\makebox(0,0){$\cdot$}}

\put(192,48){\makebox(0,0){$\cdot$}}
\put(196,44){\makebox(0,0){$\cdot$}}
\put(200,40){\makebox(0,0){$\cdot$}}
\put(204,36){\makebox(0,0){$\cdot$}}
\put(208,32){\makebox(0,0){$\cdot$}}

\put(171,80){\makebox(0,0){$\cdot$}}
\put(172,76){\makebox(0,0){$\cdot$}}
\put(174,73){\makebox(0,0){$\cdot$}}
\put(176,71){\makebox(0,0){$\cdot$}}
\put(180,70){\makebox(0,0){$\cdot$}}
\put(184,70){\makebox(0,0){$\cdot$}}
\put(188,70){\makebox(0,0){$\cdot$}}

\put(190,100){\oval(20,20)[bl]}
\put(190,90){\line(0,-1){20}}
\put(190,80){\oval(40,20)[br]}
\put(210,80){\line(0,-1){50}}
\put(210,30){\line(1,-1){10}}

\put(212,84){\makebox(0,0){$\cdot$}}
\put(214,88){\makebox(0,0){$\cdot$}}
\put(216,92){\makebox(0,0){$\cdot$}}
\put(218,96){\makebox(0,0){$\cdot$}}

\put(192,90){\makebox(0,0){$\cdot$}}
\put(194,91){\makebox(0,0){$\cdot$}}
\put(196,93){\makebox(0,0){$\cdot$}}
\put(198,96){\makebox(0,0){$\cdot$}}

\end{picture}
\end{center}

The term $t'$ in Jones normal form such that $\delta (t')$ is equivalent
to the 11-diagram whose threads are dotted is
$$c^6 h^{[3, 1]} h^{[4, 4]}h^{[7, 7]}h^{[9, 8]}h^{[10, 9]}.$$

\noindent Each solid staircase in the picture
corresponds to a block of the normal form.

However, the most simple procedure to obtain out of $t$ the term $t'$
in Jones normal form such that $t=t'$ in ${\cal K}_{n}$ 
is to draw $\delta (t)$ and recognize in it the top and
bottom slope points, from which we immediately obtain $t'$.
And to check whether $t=u$ in ${\cal K}_{n}$, it is enough
to check whether $\delta (t) \simeq \delta (u)$, which we can do without
mentioning the Jones normal form, though we relied essentially on this
normal form in order to demonstrate that $\delta '$ is an isomorphism.

\end{document}